\documentclass[11pt,a4paper]{article}

\usepackage{amsmath,amssymb,amsfonts,amsthm,bm}
\usepackage{graphicx,indentfirst,float}
\usepackage{subcaption}
\usepackage{tikz,enumerate}
\usepackage[shortlabels]{enumitem}
\usepackage[margin=2.5cm]{geometry}
\usepackage{xcolor}
\usepackage[colorlinks=true,anchorcolor=blue,filecolor=blue,
  linkcolor=blue,urlcolor=blue,citecolor=red]{hyperref}

\usepackage{tikz}
\usetikzlibrary{calc}

\newtheorem{theorem}{Theorem}[section]
\newtheorem{lemma}[theorem]{Lemma}

\newtheorem{proposition}[theorem]{Proposition}

\newtheorem{claim}{Claim}[section]
\newtheorem{definition}{Definition}[section]

\title{Tur\'{a}n number of the odd-ballooning of odd prisms\footnote{Supported by the National Natural Science Foundation of China (Nos.\,12501471, 12601692 and 12526589) and the Zhejiang Provincial Natural Science Foundation of China (No. LQN26A010013).}}
\author{{\bf Wenwen Chen$^{a,b}$},~~{\bf Longfei Fang$^{a,b}$}\thanks{Corresponding author: \texttt{lffang@chzu.edu.cn} (L. Fang)},~~{\bf Yuke Zhang$^{c}$}\\
\small $^{a}$ School of Mathematics and Finance, Chuzhou University\\
\small Chuzhou, Anhui, 239012, China\\
\small $^{b}$ Applied Mathematics Research Center, Chuzhou University,\\
\small Chuzhou, Anhui, 239012, China\\
\small $^{c}$ School of Mathematical Sciences, Zhejiang University of Technology,\\
\small Hangzhou, Zhejiang, 310032, China\\
}
\date{}
\begin{document}
\maketitle

\begin{abstract}
The Tur\'{a}n problem for odd-balloonings has been studied for various classes of skeleton graphs, including stars, paths, cycles, trees, bipartite graphs, and certain 3-chromatic graphs. The Tur\'{a}n numbers and extremal graphs for odd prisms have also been determined. Motivated by these developments, we investigate the Tur\'{a}n problem for prescribed odd-balloonings of odd prisms.
Let $k\geq1$ be fixed, and let $F^o$ be a prescribed odd-ballooning of the odd prism $C_{2k+1}\square K_2$, in which each edge is replaced by a cycle of prescribed odd length at least five. For all sufficiently large $n$, we prove that
$${\rm ex}(n,F^o)=\Big\lfloor\frac{n^2}{4}\Big\rfloor+(k+1)n-3(k+1)^2.$$
Furthermore, the unique extremal graph is
$K_{\,2k+2,\,\lfloor(n-2k-2)/2\rfloor,\,\lceil(n-2k-2)/2\rceil}$.
The proof combines the decomposition-family method and Simonovits' structural theorem with local embedding arguments to establish the extremal structure and its uniqueness.
\end{abstract}

\begin{flushleft}
\end{flushleft}

\section{Introduction}\label{sec:introduction}

Let $\mathcal{F}$ be a family of graphs. A graph is called \emph{$\mathcal{F}$-free}
if it contains no subgraph isomorphic to any member of $\mathcal{F}$. The Tur\'{a}n
number of $\mathcal{F}$, denoted by ${\rm ex}(n,\mathcal{F})$,
is the maximum number
of edges in an $n$-vertex $\mathcal{F}$-free graph.
Let ${\rm EX}(n,\mathcal{F})$
denote the family of all $n$-vertex $\mathcal{F}$-free graphs with
${\rm ex}(n,\mathcal{F})$ edges. When $\mathcal{F}=\{F\}$, we write
${\rm ex}(n,F)$ and ${\rm EX}(n,F)$.
For a graph $G$, let $V(G)$ and $E(G)$ denote its
vertex and edge sets, respectively, and write
$|G|=|V(G)|$ and $e(G)=|E(G)|$.
We write $K_{a,b,c}$ for the complete tripartite graph with part sizes $a,b,c$.

The \emph{odd-ballooning} of a graph $F$, denoted by $F^{o}$,
is obtained by replacing every edge $uv\in E(F)$ by an odd cycle
\emph{containing} the skeleton edge $uv$. Equivalently, one retains
$uv$ and adds an internally vertex-disjoint even-length $u$--$v$
path, with the new vertices for distinct edges being disjoint.
Throughout this paper, each substituted odd cycle has a prescribed
length at least five; distinct edges may have different prescribed
lengths.
Tur\'{a}n problems for odd-balloonings have been studied for
various classes of skeleton graphs. Previous results include
stars \cite{Hou2016,Hou2018,Yuan2018}, paths and cycles
\cite{Zhu2020}, and trees under additional assumptions
\cite{Zhu2023}. More recently, Geng and Xu \cite{GengXu2026}
studied odd-balloonings of lollipop graphs.

Several studies have investigated odd-ballooning of bipartite graphs.%, including complete bipartite graphs and general bipartite graphs.
The case $K_{2,2}=C_4$ is covered by the cycle results
\cite{Zhu2020}. Yan \cite{Yan2021} determined the case in which
every edge of $K_{2,3}$ is replaced by a five-cycle.
Peng and Xia \cite{Peng-2024} treated $K_{s,t}$ with
$2\le s\le t$ and $(s,t)\notin\{(2,2),(2,3)\}$,
when every replacing cycle has length at least five.
Zhai and Yuan \cite{ZhaiYuan2023} applied progressive induction
to obtain general upper and lower bounds for the odd-ballooning of
bipartite graphs.
In a different direction,
Fang, Huang, Lin and Shu \cite{FangHuangLinShu2025} determined
 the Tur\'{a}n number for the odd-ballooning of graphs
obtained by adding a universal vertex to a graph $B$,
where every component of $B$ is a non-trivial tree
or an even cycle.

Let $\square$ denote the Cartesian product. Throughout, fix an integer $k\geq 1$ and let $F=C_{2k+1}\square K_2$ be the odd prism.
He, Li and Feng \cite{HeLiFeng2025} determined the Tur\'{a}n
number and extremal graphs for the odd prism itself.
The odd prism $C_{2k+1}\square K_2$ is a natural family of non-bipartite cubic graphs consisting of two disjoint odd cycles joined by a perfect matching. Unlike the skeleton graphs treated in several existing odd-ballooning results, odd prisms are neither bipartite nor equipped with a universal vertex. These structural features make their odd-balloonings a natural subject for further investigation.
In this paper, we establish the exact Tur\'{a}n number and characterize the unique extremal graph for every prescribed odd-ballooning of the odd prism  $F$, in which
each edge $e\in E(F)$ is replaced by an odd cycle of prescribed length
$\ell_e\geq 5$.

\begin{theorem}\label{thm:main}
Let $k\geq1$ be an integer and let $F=C_{2k+1}\square K_2$.
Let $F^o$ be a prescribed odd-ballooning of $F$ in which each edge
$e\in E(F)$ is replaced by a cycle of prescribed odd length
$\ell_e\geq5$. Then, for all sufficiently large $n$, we have
\[
{\rm ex}(n,F^o)=\Big\lfloor\frac{n^2}{4}\Big\rfloor
 +(k+1)n-3(k+1)^2.
\]
Moreover, the unique extremal graph is
$K_{\,2k+2,\,\lfloor(n-2k-2)/2\rfloor,\,\lceil(n-2k-2)/2\rceil}$.
\end{theorem}

The remainder of the paper is organized as follows. In Section~\ref{sec:decomposition}, we establish the relevant coloring and decomposition properties and recall the structural theorem needed later. In Section~\ref{sec:main-proof}, we use these ingredients to complete the proof of the main extremal result.

\section{Decomposition family and Simonovits' theorem}\label{sec:decomposition}

Write $E_t=\overline{K_t}$ for the
empty graph on $t$ vertices, and let $G\vee J$ and $G\cup J$ denote
the join and disjoint union, respectively.
Define $T_{n,t}$ to be the balanced complete
$t$-partite graph on $n$ vertices, whose part sizes
differ by at most one.
The decomposition-family method goes back to Simonovits
\cite{Simonovits1974,Simonovits1983};
we use the following version due to Liu \cite{Liu2013}.

\begin{definition}\label{definition2.1A}
Let $H$ be a graph with $\chi(H)=p+1\ge 3$.
The \textbf{decomposition family} $\mathcal{M}(H)$ consists of the
minimal graphs $M$ for which there exists a constant $t=t(H)$ such that
\[
H\subseteq (M\cup E_t)\vee T_{(p-1)t,p-1}.
\]
\end{definition}
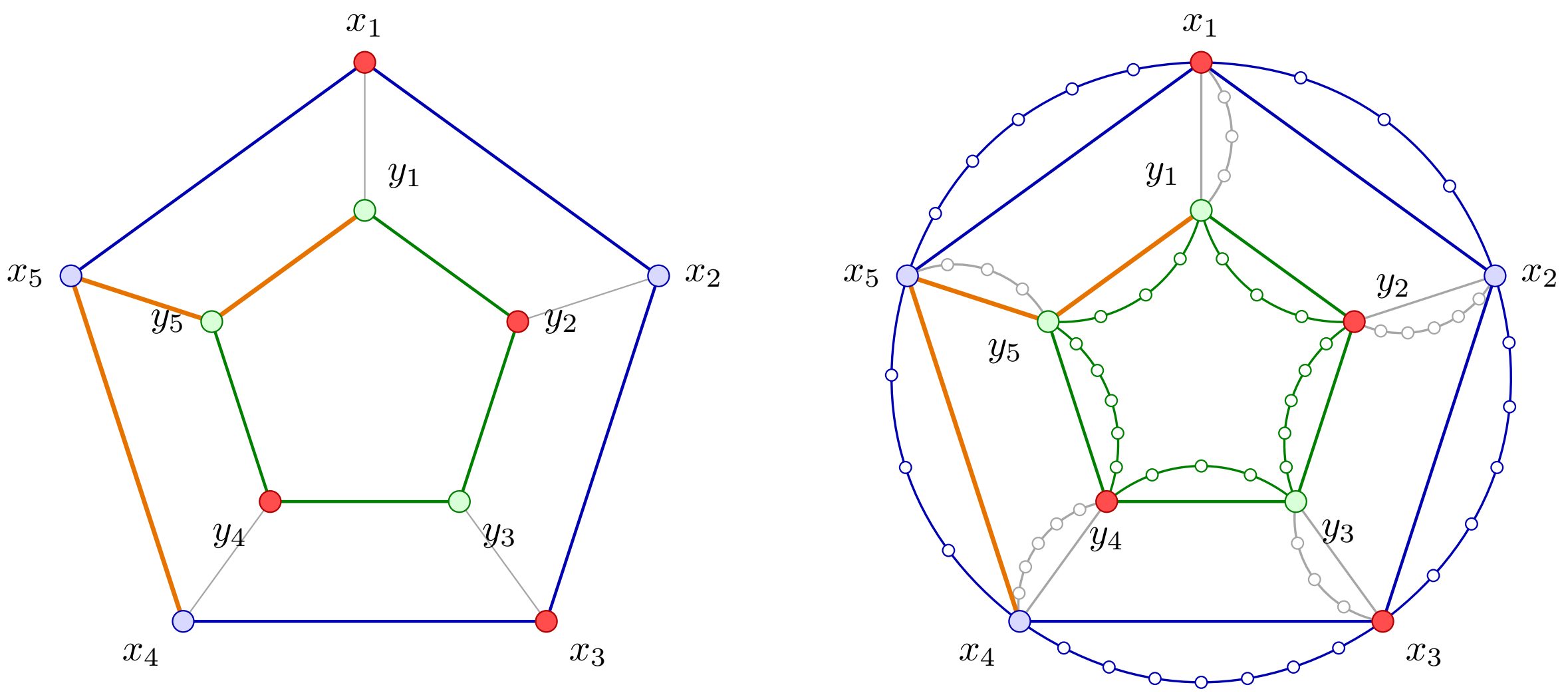
\begin{figure}[htbp]
    \centering
    \begin{subfigure}[t]{0.48\linewidth}
        \centering
        \resizebox{\linewidth}{!}{%
            \begin{tikzpicture}[
            scale=1.2,
            vertex/.style={circle, draw=black, fill=white, inner sep=2pt},
            outervertex/.style={circle, draw=blue!70!black, fill=blue!15, inner sep=2pt},
            innervertex/.style={circle, draw=green!50!black, fill=green!15, inner sep=2pt},
            Iouter/.style={circle, draw=red!70!black, fill=red!70, inner sep=2pt},
            Iinner/.style={circle, draw=red!70!black, fill=red!70, inner sep=2pt},
            every node/.style={font=\small}
            ]

            % Use identical drawing bounds to keep both graphs at the same scale.
            \path[use as bounding box] (-3,-2.65) rectangle (3,2.90);
            % Outer cycle: x_1,\dots,x_5
            \coordinate (x1) at (90:2.4);
            \coordinate (x2) at (18:2.4);
            \coordinate (x3) at (-54:2.4);
            \coordinate (x4) at (-126:2.4);
            \coordinate (x5) at (162:2.4);

            % Inner cycle: y_1,\dots,y_5
            \coordinate (y1) at (90:1.25);
            \coordinate (y2) at (18:1.25);
            \coordinate (y3) at (-54:1.25);
            \coordinate (y4) at (-126:1.25);
            \coordinate (y5) at (162:1.25);

            % Outer 5-cycle
            \draw[blue!70!black, thick] (x1)--(x2)--(x3)--(x4)--(x5)--cycle;

            % Inner 5-cycle
            \draw[green!50!black, thick] (y1)--(y2)--(y3)--(y4)--(y5)--cycle;

            % Rungs
            \foreach \i in {1,...,5}
                \draw[gray!70] (x\i)--(y\i);

            % Highlight the P_4 in F-I: x_4 x_5 y_5 y_1
            \draw[orange!90!black, very thick]
                (x4)--(x5)--(y5)--(y1);

            % Vertices in I = {x_1,x_3,y_2,y_4}
            \node[Iouter, label=above:$x_1$] at (x1) {};
            \node[outervertex, label=right:$x_2$] at (x2) {};
            \node[Iouter, label=below right:$x_3$] at (x3) {};
            \node[outervertex, label=below left:$x_4$] at (x4) {};
            \node[outervertex, label=left:$x_5$] at (x5) {};

            \node[innervertex, label=above right:$y_1$] at (y1) {};
            \node[Iinner, label=right:$y_2$] at (y2) {};
            \node[innervertex, label=below right:$y_3$] at (y3) {};
            \node[Iinner, label=below left:$y_4$] at (y4) {};
            \node[innervertex, label=left:$y_5$] at (y5) {};

            \end{tikzpicture}%
        }
        \caption{The prism $F=C_5\square K_2$ ($k=2$).}
        \label{fig:1}
    \end{subfigure}\hfill
    \begin{subfigure}[t]{0.48\linewidth}
        \centering
        \resizebox{\linewidth}{!}{%
            \begin{tikzpicture}[
                scale=1.2,
                outervertex/.style={circle, draw=blue!70!black,
                    fill=blue!15, inner sep=2pt},
                innervertex/.style={circle, draw=green!50!black,
                    fill=green!15, inner sep=2pt},
                Ivertex/.style={circle, draw=red!70!black,
                    fill=red!70, inner sep=2pt},
                newvertex/.style={circle, fill=white, inner sep=1.1pt},
                every node/.style={font=\small}
            ]

            % Use identical drawing bounds to keep both graphs at the same scale.
            \path[use as bounding box] (-3,-2.65) rectangle (3,2.90);

            % Arguments: u, v, prescribed odd cycle length, bend angle, color.
            % Keep uv and add a path with ell-1 edges and ell-2 new vertices.
            \newcommand{\AddBalloon}[5]{%
                \draw[#5, semithick] (#1) to[bend left=#4] (#2);
                \pgfmathtruncatemacro{\nnew}{#3-2}
                \foreach \s in {1,...,\nnew}{
                    \pgfmathsetmacro{\t}{\s/(#3-1)}
                    \path (#1) to[bend left=#4]
                        node[pos=\t, newvertex, draw=#5] {} (#2);
                }
            }

            % The ten skeleton vertices, in the same layout as the original F.
            \foreach \i in {1,...,5}{
                \pgfmathsetmacro{\ang}{90-72*(\i-1)}
                \coordinate (x\i) at (\ang:2.4);
                \coordinate (y\i) at (\ang:1.25);
            }

            % Retain ALL fifteen skeleton edges.
            \draw[blue!70!black, thick]
                (x1)--(x2)--(x3)--(x4)--(x5)--cycle;
            \draw[green!50!black, thick]
                (y1)--(y2)--(y3)--(y4)--(y5)--cycle;
            \foreach \i in {1,...,5}
                \draw[gray!70, semithick] (x\i)--(y\i);

            % Prescribed lengths ell_{x_i x_j}: edit the third entry in each triple.
            % The added outer paths bend away from the outer pentagon.
            \foreach \i/\j/\L in {1/2/5, 2/3/7, 3/4/9, 4/5/5, 5/1/7}{
                \AddBalloon{x\i}{x\j}{\L}{35}{blue!70!black}
            }

            % Prescribed lengths ell_{y_i y_j}.
            % The added inner paths bend into the inner pentagon.
            \foreach \i/\j/\L in {1/2/5, 2/3/7, 3/4/5, 4/5/7, 5/1/5}{
                \AddBalloon{y\i}{y\j}{\L}{-40}{green!50!black}
            }

            % Prescribed lengths ell_{x_i y_i}.
            % These paths lie in the annulus between the two pentagons.
            \foreach \i/\L in {1/5, 2/7, 3/5, 4/7, 5/5}{
                \AddBalloon{x\i}{y\i}{\L}{45}{gray!70}
            }

            % Preserve the highlighted skeleton path from the original figure.
            \draw[orange!90!black, very thick]
                (x4)--(x5)--(y5)--(y1);

            % Draw the skeleton vertices last.
            \node[Ivertex, label=above:\(x_1\)] at (x1) {};
            \node[outervertex, label=right:\(x_2\)] at (x2) {};
            \node[Ivertex, label=below right:\(x_3\)] at (x3) {};
            \node[outervertex, label=below left:\(x_4\)] at (x4) {};
            \node[outervertex, label=left:\(x_5\)] at (x5) {};

            % Place inner labels to the side of the rungs.
            \node[innervertex, label=125:\(y_1\)] at (y1) {};
            \node[Ivertex, label=53:\(y_2\)] at (y2) {};
            \node[innervertex, label=-19:\(y_3\)] at (y3) {};
            \node[Ivertex, label=-91:\(y_4\)] at (y4) {};
            \node[innervertex, label=197:\(y_5\)] at (y5) {};

            \end{tikzpicture}%
        }
        \caption{A prescribed odd-ballooning $F^o$.}
        \label{fig:odd-ballooning-C5-prism}
    \end{subfigure}

    \caption{The prism $F$ and its prescribed odd-ballooning $F^o$.}
    \label{fig:C5-prism-and-ballooning}
\end{figure}

To prepare for the analysis of $\mathcal M(F^o)$, we first record
several elementary properties of the odd prism $F$ that will also
be useful in constructing the lower-bound example.
Label the vertices of $F$ as $x_i$ and $y_i$, for $1\le i\le 2k+1$, with all indices taken modulo $2k+1$. Its edges consist of the two rail cycles together with the rungs $x_iy_i$. Figure~\ref{fig:1} illustrates this structure and the labeling in the case $k=2$, and Figure~\ref{fig:odd-ballooning-C5-prism} shows its prescribed odd-ballooning $F^o$ where the added paths have pairwise disjoint sets of internal vertices, shown as small unlabelled circles. %Thus $|F|=4k+2$ and $e(F)=6k+3$.

\begin{lemma}\label{lem:2.1}
The graph $F$ satisfies the following properties:

{\rm (i)}  $\alpha(F)\le 2k$, where $\alpha(F)$ denotes the independence number of $F$;

{\rm (ii)} there exists an independent set $I$ of size $2k$ such that
$F-I\cong P_4\cup(2k-2)K_1$, where the second term is empty when $k=1$;

{\rm (iii)}
there exists  a partition $V(F)=S\mathbin{\dot\cup}T$ with
$|S|=|T|=2k+1$ such that each of $F[S]$ and $F[T]$ has exactly one edge.
\end{lemma}

\begin{proof}
{\rm (i)}
An independent set meets each rung in at most one vertex.
Suppose that an independent set contains exactly one vertex from each of the $2k+1$ rungs.
If the chosen vertex on one rung belongs to one rail cycle, then the chosen vertex on the next rung must belong to the other rail cycle, since vertices on consecutive rungs belonging to the same rail cycle are adjacent.
Thus the choices are forced to alternate between the two rail cycles.
Since there are $2k+1$ rungs, the chosen vertices on the first
and the last rungs must lie on the same rail cycle. These two vertices
occupy consecutive positions on that rail cycle and are therefore
adjacent, which contradicts the independence of the set.
Consequently, $\alpha(F)\le 2k$.

\vspace{2mm}
{\rm (ii)}
Set
$I=\{x_1,x_3,\ldots,x_{2k-1}\} \cup\{y_2,y_4,\ldots,y_{2k}\}$.
This set is independent and has size $2k$. The remaining edges are
exactly the three edges of the path
$x_{2k}\,x_{2k+1}\,y_{2k+1}\,y_1$.
Thus, $F-I\cong P_4\cup(2k-2)K_1$.

\vspace{2mm}
{\rm (iii)}
Set $S=I\cup\{x_{2k+1}\}$ and $T=V(F)\setminus S$.
The only edge of $F[S]$ is $x_{2k+1}x_1$, and that of $F[T]$ is $y_{2k+1}y_1$.
\end{proof}

We now use these properties of $F$ to determine the relevant
coloring parameters of $F^o$ and to obtain the lower bound in
Theorem~\ref{thm:main}.
For a three-chromatic graph $L$, let $\sigma(L)$ be the minimum size
of a color class over all its proper three-colorings.
For convenience, set
\[
g_k(n):=
\Big\lfloor\frac{n^2}{4}\Big\rfloor
+(k+1)n-3(k+1)^2.
\]

\begin{lemma}\label{lem:2.2}
The graph $F^o$ satisfies the following properties:

{\rm (i)} we have $\chi(F^o)=3$;

{\rm (ii)}  $\sigma(F^o)=2k+3$;

{\rm (iii)} $E_{2k+2}\vee T_{n-2k-2,2}$ is
$F^o$-free whenever it is defined. In particular ${\rm ex}(n,F^o)\ge g_k(n)$.
\end{lemma}
\begin{proof}

{\rm (i)} Since every edge of $F$ is replaced by an odd cycle, the graph $F^o$ contains an odd cycle. Hence $\chi(F^o)\ge 3$.
It remains to show that $F^o$ is three-colorable. Let
\[
I=\{x_1,x_3,\ldots,x_{2k-1}\}
 \cup\{y_2,y_4,\ldots,y_{2k}\}.
\]
Since $F-I$ is bipartite, we may properly color $F-I$ with colors one and two,
and assign color three to every vertex of $I$.

For each edge $e=uv\in E(F)$, let
$u z_1z_2\cdots z_{\ell_e-2}v$
be the additional $u$--$v$ path which, together with the edge $uv$, forms the odd cycle replacing $e$. Since $\ell_e$ is odd, $\ell_e-2$ is odd.
Moreover, since $I$ is an independent set of $F$ and $uv\in E(F)$, at most one of $u$ and $v$ belongs to $I$.

If exactly one of $u$ and $v$ belongs to $I$, say $u\in I$,
then $v$ is assigned color one or two.
Color $z_1,\ldots,z_{\ell_e-2}$ alternately with colors one and two, starting with the color different from that of $v$. Since $\ell_e-2$ is odd, the last vertex $z_{\ell_e-2}$ also receives the color different from that of $v$. Thus the whole additional path is properly colored.

If $e=uv$ is an edge of $F-I$, then $u$ and $v$ receive different colors from $\{1,2\}$.
Color $z_1,\ldots,z_{\ell_e-3}$ alternately with colors one and two,
starting with the color different from that of $u$,
and assign color three to $z_{\ell_e-2}$.
Again, the resulting coloring of the additional path is proper.

Since the additional paths corresponding to distinct edges of $F$ are internally vertex-disjoint,
these colorings are compatible and together give a proper three-coloring of $F^o$.
Therefore $\chi(F^o)\le 3$.
Combining this with $\chi(F^o)\ge 3$,
we conclude that $\chi(F^o)=3$.

\vspace{2mm}
{\rm (ii)}
We first prove that $\sigma(F^o)\le 2k+3$.
Count the number of vertices receiving color three in the three-coloring constructed above. First, every vertex of $I$ is assigned color three, so these vertices contribute exactly $|I|=2k$ vertices to the third color class.

Next, consider the new vertices introduced on the replacing paths.
If a skeleton edge $e\in E(F)$ is incident with a vertex of $I$, then all new internal vertices on the replacing path corresponding to $e$ are colored only with colors one and two. Thus such an edge contributes no new vertex of color three.

On the other hand, if $e\in E(F-I)$, then both endpoints of $e$ lie outside $I$. In this case, by our coloring rule, exactly one new internal vertex on the replacing path corresponding to $e$ is assigned color three, while all the remaining new vertices on that path receive colors one and two. Since the replacing paths corresponding to distinct edges of $F$ are internally vertex-disjoint, these newly colored vertices are all distinct. Therefore, the number of new vertices receiving color three is exactly $e(F-I)$.

By Lemma~\ref{lem:2.1}(ii), we have $F-I\cong P_4\cup(2k-2)K_1$. Hence $e(F-I)=e(P_4)=3$. Consequently, the third color class has size
$|I|+e(F-I)=2k+3$.
Therefore, the three-coloring constructed above has a color class of size exactly $2k+3$, and hence $\sigma(F^o)\le 2k+3$.

 To establish the reverse inequality, let $Z$ be an arbitrary color class in an arbitrary proper three-coloring of $F^o$. It is enough to show that $|Z|\ge 2k+3$. This will imply $\sigma(F^o)\ge 2k+3$, and hence $\sigma(F^o)=2k+3$.

Set $U=Z\cap V(F)$.
Since vertices in the same color class are pairwise nonadjacent,
$U$ is an independent set in $F$.
By Lemma \ref{lem:2.1}(i), we have
$|U|\le \alpha(F)\le 2k$.

After deleting the color class $Z$, the graph $F^o-Z$ is colored with the remaining two colors and is therefore bipartite. Consider an edge $e=uv\in E(F-U)$. Since neither $u$ nor $v$ belongs to $U$, neither endpoint of $e$ belongs to $Z$. Recall that in $F^o$, the edge $e$ lies on the odd cycle that replaces $e$. If none of the new internal vertices of this replacing cycle belonged to $Z$, then the entire odd cycle would remain in $F^o-Z$, contradicting the fact that $F^o-Z$ is bipartite. Hence, for every edge $e\in E(F-U)$, the corresponding replacing cycle contains at least one new vertex of $Z$.
The sets of new vertices belonging to distinct replacing cycles are pairwise disjoint. Therefore, in addition to the $|U|$ original vertices of $F$ contained in $Z$, the color class $Z$ contains at least one new vertex for each edge of $F-U$. Thus $|Z|\ge |U|+e(F-U)$.

Now $F=C_{2k+1}\square K_2$ is cubic and has $e(F)=3(2k+1)=6k+3$. Since $U$ is independent, no edge of $F$ has both endpoints in $U$. Hence deleting the vertices of $U$ removes exactly $3|U|$ edges, and consequently $e(F-U)=6k+3-3|U|$. It follows that
$$|Z|\ge |U|+e(F-U)=6k+3-2|U|\ge 6k+3-4k=2k+3,$$
 as desired.

\vspace{2mm}
{\rm (iii)}
Suppose to the contrary that there exists an embedding of $F^o$ into $E_{2k+2}\vee T_{n-2k-2,2}$. Let $V_0$ denote the vertex set of $E_{2k+2}$, and let $V_1$ and $V_2$ be the two partite sets of $T_{n-2k-2,2}$. By construction, each of $V_0,V_1,V_2$ is an independent set in the host graph.

Pull these three parts back through the embedding of $F^o$. Namely, assign color $i$ to a vertex of $F^o$ whenever its image lies in $V_i$, for $i=0,1,2$. Since each $V_i$ is independent, two adjacent vertices of $F^o$ cannot be mapped into the same $V_i$. Hence this assignment gives a proper three-coloring of $F^o$.

Consider the color class corresponding to $V_0=V(E_{2k+2})$. Since $V_0$ contains exactly $2k+2$ vertices and the embedding is injective, at most $2k+2$ vertices of $F^o$ can be mapped into $V_0$. Thus this proper three-coloring of $F^o$ has a color class of size at most $2k+2$.

However, by the definition of $\sigma(F^o)$
and the equality $\sigma(F^o)=2k+3$ established above,
every color class in every proper three-coloring of $F^o$ has size at least $2k+3$.
This contradicts the existence of a color class of size at most $2k+2$.
Therefore $F^o$ cannot be embedded into $E_{2k+2}\vee T_{n-2k-2,2}$,
and hence $E_{2k+2}\vee T_{n-2k-2,2}$ is $F^o$-free.

It remains to compute its number of edges. Since $2k+2$ is even, we have
\begin{align*}
e(E_{2k+2}\vee T_{n-2k-2,2})
&=(2k+2)(n-2k-2)
  +\Big\lfloor\frac{(n-2k-2)^2}{4}\Big\rfloor\\
&=\Big\lfloor\frac{n^2}{4}\Big\rfloor
  +(k+1)n-3(k+1)^2=g_k(n).
\end{align*}
Indeed,
$\frac{(n-2k-2)^2}{4}
=\frac{n^2}{4}-(k+1)n+(k+1)^2$,
and the last two terms are integers. Therefore $E_{2k+2}\vee T_{n-2k-2,2}$ is an $F^o$-free graph with $g_k(n)$ edges, which proves that ${\rm ex}(n,F^o)\ge g_k(n)$.
\end{proof}

\subsection{Decomposition family and its construction}

Since $\chi(F^o)=3$ by Lemma~\ref{lem:2.2}(i), we have $p=2$ in Definition~\ref{definition2.1A}.
Hence, for $F^o$, the definition specializes to $F^o\subseteq(M\cup E_t)\vee E_t$.
We next recall the cracking characterization that will be used throughout the proof.
Let $F$ be a graph, and let $u\in V(F)$ be a vertex of degree
$d$ with $N_F(u)=\{v_1,\ldots,v_d\}$. Cracking $u$ replaces it
by $d$ distinct new vertices $u_1,\ldots,u_d$. For each
$i\in\{1,2,\dots,d\}$, the edge $uv_i$ is replaced independently in one
of the following two ways.

\textbf{Type I:} replace $uv_i$ by $u_iv_i$; that is, the new vertex
$u_i$ is joined to the original neighbor $v_i$.

\textbf{Type II:} replace $uv_i$ by $u_iw_i$, where $w_i$ is a new
vertex; that is, the original neighbor $v_i$ is detached from this
edge and a new isolated edge is created between $u_i$ and $w_i$.

For an independent set $U\subseteq V(F)$, crack all vertices in
$U$ simultaneously and delete isolated vertices. Denote the
family of all resulting graphs by $\mathcal C(F,U)$, and let
\[
\mathcal C(F)=\bigcup_{\substack{U\subseteq V(F)\\U\text{ independent}}}
\mathcal C(F,U).
\]

Let
$V(F)=\{x_1,\ldots,x_{2k+1},y_1,\ldots,y_{2k+1}\}$
and
\[
E(F)=\{x_ix_{i+1},y_iy_{i+1}:1\le i\le 2k+1\}
\cup\{x_iy_i:1\le i\le 2k+1\}.
\]
The edges $x_iy_i$ ($1\le i\le 2k+1$) are called
the \emph{rung edges} of $F$.
Thus $e(F)=6k+3$. For an independent set $U\subseteq V(F)$, let
\[
\delta_{F}(U)=\{uv\in E(F):u\in U,\ v\notin U\}.
\]
Since $U$ is independent, every edge in $\delta_{F}(U)$ has exactly one endpoint
in $U$. For each map
$\tau:\delta_{F}(U)\longrightarrow\{\mathrm{I},\mathrm{II}\}$,
define $F(U,\tau)$ as follows.
For every edge $e=uv\in\delta_{F}(U)$, where
$u\in U$ and $v\notin U$, replace $u$ on $e$ by a new vertex $u_e$. If
$\tau(e)=\mathrm{I}$, replace $e$ by $u_ev$; if
$\tau(e)=\mathrm{II}$, replace $e$ by $u_ew_e$, where $w_e$ is a new
vertex. All vertices introduced for distinct edges are distinct, every edge of
$F-U$ is left unchanged, and isolated vertices are deleted at the end.

We quote the cracking lemma of Fang, Huang, Lin and Shu
\cite[Lemma~2.1]{FangHuangLinShu2025}.

\begin{lemma}[{\cite[Lemma~2.1]{FangHuangLinShu2025}}]\label{lem:2.3}
Let $F$ be a graph with $\chi(F)\in\{2,3\}$.
Then $\mathcal{C}(F)=\mathcal{M}(F^{o})$.
\end{lemma}

We now apply Lemma~\ref{lem:2.3} to the skeleton $F$.

\begin{proposition}\label{prop:2.4}
Let $F^o$ be the prescribed odd-ballooning of $F=C_{2k+1}\square K_2$. Then
\[
\mathcal M(F^o)=
\bigcup_{\substack{U\subseteq V(F)\\U\text{ independent}}}
\bigl\{F(U,\tau)~|~\tau:\delta_F(U)\to\{\mathrm I,\mathrm {II}\}\bigr\}.
\]
In particular, every member of $\mathcal{M}(F^o)$ is obtained by choosing
an independent set $U$ of $F$, cracking every vertex of $U$, and independently choosing a Type I or Type II replacement for each edge incident with $U$.
\end{proposition}

\begin{proof}
Since $F$ contains $C_{2k+1}$ and $F\subseteq F^o$,
Lemma~\ref{lem:2.2}(i) gives
$3\le\chi(F)\le\chi(F^o)=3$, which implies $\chi(F)=3$.
Thus Lemma~\ref{lem:2.3} applies and gives
$\mathcal M(F^o)=\mathcal C(F)$.

Now fix an independent set $U\subseteq V(F)$. Since $U$ is independent,
no edge of $F$ has both endpoints in $U$. Therefore, every edge of $F$
either has both endpoints in $V(F)\setminus U$ or has exactly one endpoint
in $U$. The former edges remain unchanged when the vertices of $U$ are
cracked, whereas the latter edges are precisely the edges in
$\delta_{F}(U)$.

For each edge $e\in\delta_{F}(U)$, the cracking operation replaces $e$ in
exactly one of two possible ways, namely by a Type~I or a Type~II replacement.
Moreover, these choices can be made independently for the different edges
in $\delta_{F}(U)$. Thus all choices involved in simultaneously cracking
the vertices of $U$ are encoded by a map
$\tau:\delta_{F}(U)\to\{\mathrm{I},\mathrm{II}\}$, where $\tau(e)$
records the type of replacement chosen for $e$.

Once $U$ and $\tau$ are fixed, the resulting graph is uniquely
determined up to isomorphism:
the edges with both endpoints outside $U$ remain unchanged, while each edge
$e\in\delta_{F}(U)$ is replaced according to $\tau(e)$, with distinct
new vertices used for distinct skeleton edges. After deleting any isolated
vertices created in the process, the resulting graph is precisely
$F(U,\tau)$.

Conversely, every map
$\tau:\delta_{F}(U)\to\{\mathrm{I},\mathrm{II}\}$ specifies one of the
two allowed replacements for each edge incident with $U$, and hence
determines a valid simultaneous cracking of all vertices in $U$. Therefore,
$\mathcal{C}(F,U)=
\{F(U,\tau)~|~
\tau:\delta_{F}(U)\to\{\mathrm{I},\mathrm{II}\}\}$.

Finally, taking the union over all independent sets
$U\subseteq V(F)$ gives the asserted description of $\mathcal{C}(F)$.
Together with $\mathcal M(F^o)=\mathcal C(F)$ from
Lemma~\ref{lem:2.3}, this proves the proposition.
\end{proof}

The full description of the decomposition family immediately yields
the linear forest $P_4\cup 6kK_2$ used later in the proof,
so no separate minimality argument is needed.
For convenience, let
\[
f:=|F^o|=4k+2+\sum_{e\in E(F)}(\ell_e-2).
\]

\begin{lemma}\label{lem:2.5}
We have $P_4\cup6kK_2\in\mathcal{M}(F^o)$. In particular, $\mathcal{M}(F^o)$ contains a linear
forest, and
\begin{equation}\label{eq:2A}
F^o\subseteq(P_4\cup 6kK_2 \cup E_f)\vee E_f.
\end{equation}
\end{lemma}
\begin{proof}
Take the independent set $I$ from Lemma~\ref{lem:2.1}. Since $F$ is
cubic and $|I|=2k$, exactly $6k$ edges of $F$ have one endpoint in $I$.
By Lemma~\ref{lem:2.1}, after deleting $I$ the only nontrivial
component is a copy of $P_4$.

Crack every vertex of $I$ and choose Type~II for each of the $6k$ edges
in $\delta_F(I)$. Each such edge becomes an isolated edge with two fresh
endpoints, while the three edges of $F-I$ remain the $P_4$. After deleting
isolated vertices, the resulting graph is exactly
$P_4\cup 6kK_2$.
Proposition~\ref{prop:2.4} therefore gives
$P_4\cup 6kK_2\in\mathcal{M}(F^o)$.

We use the following immediate consequence of the definition of the decomposition family. Since $\chi(F^o)=3$, if $M\in\mathcal M(F^o)$, then there exists an integer $t$ such that
$F^o\subseteq (M\cup E_t)\vee E_t$.
Recall that $f=|F^o|$. In any such embedding, at most $f$ vertices of each copy of $E_t$ are used, since $F^o$ has only $f$ vertices. Hence the used vertices in each independent part can be mapped injectively into a copy of $E_f$. Therefore
$F^o\subseteq (M\cup E_f)\vee E_f$.
Taking $M=P_4\cup6kK_2$, which belongs to $\mathcal M(F^o)$, we obtain
$F^o\subseteq(P_4\cup6kK_2\cup E_f)\vee E_f$,
as required.
\end{proof}

We first recall the following class of graphs introduced by Simonovits \cite{Simonovits1974}, which will be useful in describing the structure of extremal graphs.

\begin{definition}\label{def:2.2}\emph{(\cite{Simonovits1974})}
Denote by $\mathcal{D}(n,p,q)$ the family of $n$-vertex graphs $G$ satisfying the
following symmetry condition:

\vspace{1mm}
{\rm (i)} There exists a set $S\subseteq V(G)$ with $|S|\leq q$ such that
$G-S=G'=\vee_{i=1}^{p}G_i$,
where $\big||G_i|-\frac{n}{p}\big|\leq q$ for every $i\in \{1,2,\ldots,p\}$.

\vspace{1mm}
{\rm (ii)} For each $i\in\{1,2,\ldots,p\}$,
there exists a connected graph $H_i$ such that $G_i$
is the disjoint union of $k_i$ copies of $H_i$,
where $k_i=\frac{|G_i|}{|H_i|}$ and $|H_i|\leq q$.
Moreover, for any two copies $H_i'$ and $H_i''$ of $H_i$ in $G_i$,
there exists an isomorphism
$\phi:V(H_i')\to V(H_i'')$ such that, for every $u\in V(H_i')$ and every $v\in S$, $uv\in E(G)$ if and only if $\phi(u)v\in E(G)$.
\end{definition}

Our proof of the main result relies on a result of Simonovits \cite{Simonovits1974}.

\begin{lemma}\emph{(\cite[Theorems~1 and~2]{Simonovits1974})}\label{lem:2.6}
Let $H$ be a graph with $\chi(H)=p+1\geq 3$.
If $\mathcal{M}(H)$ contains a linear forest,
then there exist $q=q(H)$ and $n_0=n_0(q)$ such that
$\mathcal{D}(n,p,q)$ contains a graph $G\in \mathrm{EX}(n,H)$ for all $n\geq n_0$.
Furthermore, if $G$ is the only extremal graph in $\mathcal{D}(n,p,q)$,
then it is the unique graph in $\mathrm{EX}(n,H)$ for every sufficiently large $n$.
\end{lemma}

\section{Proof of the main theorem}\label{sec:main-proof}

\begin{proof}[\textbf{Proof of Theorem~\ref{thm:main}}]
All implicit constants below depend only on the fixed graph $F^o$,
and $n$ is assumed sufficiently large. By Lemmas~\ref{lem:2.2}(i),
\ref{lem:2.5}, and~\ref{lem:2.6}, there is a constant $q$
such that ${\rm EX}(n,F^o)\cap \mathcal{D}(n,2,q)\ne\varnothing$.
We determine every extremal graph in this intersection and then apply
the uniqueness assertion of Lemma~\ref{lem:2.6}.

Fix $G\in {\rm EX}(n,F^o)\cap\mathcal D(n,2,q)$. By
Definition~\ref{def:2.2}, there is a set $R_0$ such that
\begin{equation*}%\label{eq:3A}
G-R_0=G_1\vee G_2,\qquad |R_0|\le q,\qquad
\left||G_i|-\frac n2\right|\le q\quad(i=1,2),
\end{equation*}
where each $G_i$ consists of copies of a fixed connected graph $H_i$
of order at most $q$. These copies are symmetric with respect to
$R_0$ in the sense of Definition~\ref{def:2.2}(ii).
Write $r_0=|R_0|$ and $n_i=|G_i|$. Lemma~\ref{lem:2.2}(iii) gives
\begin{equation}\label{eq:4A}
e(G)\ge g_k(n)=\frac{n^2}{4}+(k+1)n+O(1).
\end{equation}

We first show that each $G_i$ is $P_4$-free.
Suppose to the contrary that $G_i$ contains a copy of $P_4$ for some $i\in\{1,2\}$.
Recall that $G_i$ is the disjoint union of copies of the connected graph
$H_i$. Hence, if $G_i$ contains a copy of $P_4$, then $H_i$ itself
contains a copy of $P_4$, and in particular, $H_i$ contains an edge.
Since $n$ is sufficiently large, $G_i$ contains at least $6k+1$
vertex-disjoint copies of $H_i$. Taking a copy of $P_4$ from one copy
of $H_i$ and one edge from each of another $6k$ copies of $H_i$, we obtain
$P_4\cup 6kK_2\subseteq G_i$.
Since each $G_j$ is a disjoint union of copies of the fixed graph $H_j$ of bounded order, the number of such copies tends to infinity with $n$.
Hence, by selecting one vertex from each of \(f\) distinct unused copies, we obtain an independent set of size \(f\) in each required part.
Therefore,
by Lemma~\ref{lem:2.5},
$$F^o\subseteq (P_4\cup 6kK_2\cup E_f)\vee E_f\subseteq G,$$
a contradiction. Thus, $G_i$ is $P_4$-free for each $i\in\{1,2\}$.
Therefore, for each $i\in\{1,2\}$,
by the Erd\H{o}s--Gallai theorem \cite{ErdosGallai1959},  we obtain
\begin{align}\label{eq:5A}
e(G_i)\le |G_i|=n_i.
\end{align}

We next record two local embedding observations. Let $H_5$ be the
odd-ballooning of $F$ in which every replacing cycle has length five.
Since $|F|=4k+2$ and $e(F)=6k+3$, we have
\[
h_5:=|H_5|
=|F|+3e(F)
=4k+2+3(6k+3)
=22k+11.
\]
It is clear that $h_5\le f$. Set $\beta:=2f+10$.

\begin{claim}\label{claim:3.2}
Let $X,Y\subseteq V(G)$ be disjoint, completely joined sets with
$|X|,|Y|\ge f$. If a copy of $H_5$ has a non-skeleton $X$--$Y$
edge in each of its replacing five-cycles, then it can be extended
to a copy of $F^o$ with the same skeleton.
\end{claim}

\begin{proof}
For $e\in E(F)$, set $d_e=(\ell_e-5)/2$. If $d_e=0$, keep its five-cycle
unchanged. Otherwise, for its chosen edge $xy$, with $x\in X$ and $y\in Y$,
replace that edge by
\[
x\,y_{e,1}x_{e,1}\,y_{e,2}x_{e,2}\cdots
y_{e,d_e}x_{e,d_e}\,y,
\]
using fresh $x_{e,j}\in X$ and $y_{e,j}\in Y$. Every new edge exists
because $X$ and $Y$ are completely joined, and the replacement increases
the cycle length by exactly $2d_e$ without changing the skeleton edge.

Choose the new vertices outside the original copy of $H_5$ and distinctly
for different skeleton edges. The copy of $H_5$ occupies at most $h_5$
vertices in total, so each set has at least $f-h_5$ unused vertices.
Moreover,
\[
\sum_{e\in E(F)}d_e
=\frac12\sum_{e\in E(F)}(\ell_e-5)
=\frac{f-h_5}{2}.
\]
Thus exactly $(f-h_5)/2$ new vertices are required from each set.
Since each set contains at least $f-h_5$ unused vertices,
the required vertices can be chosen so that all such choices are pairwise disjoint.
The resulting cycles have the
prescribed lengths and intersect only in their prescribed skeleton
endpoints.
\end{proof}

\begin{claim}\label{claim:3.3}
Let $X,Y,C\subseteq V(G)$ be pairwise disjoint, with $|X|,|Y|\ge \beta$.
Suppose that $X$ is complete to $Y$, and that $C$ is complete to $X\cup Y$.
If $G-C$ contains a matching $\mathcal E$ of at most
three edges such that
\[
|C|+|\mathcal E|\ge2k+3,
\]
and the two endpoints of each edge of $\mathcal E$ are both complete
to $X$ or both complete to $Y$, then $F^o\subseteq G$.
The empty matching is allowed.
\end{claim}

\begin{proof}
Since $|\mathcal E|\le 3$, the size condition gives $|C|\ge 2k$.
Let $W$ be the set of all endpoints of the matching $\mathcal E$.
Then $|W|\le 6$. We reserve all vertices of $W$ throughout the
construction, using them only as endpoints of their assigned matching
edges.
Embed the independent set $I$ from Lemma~\ref{lem:2.1}(ii)
in $2k$ vertices of $C$. Embed $F-I$ in $X\cup Y$ according to
its two-coloring, avoiding $W$.
Its three edges form a copy of $P_4$,
and each of these edges has one endpoint in $X$ and the other in $Y$.

Recall that $F-I\cong P_4\cup(2k-2)K_1$, so the nontrivial component of $F-I$ is a copy of $P_4$ with exactly three skeleton edges. Label these three edges by $e_1,e_2,e_3$, and write $\mathcal E=\{f_1,\ldots,f_s\}$, where $s=|\mathcal E|\le3$. After relabeling $e_1,e_2,e_3$ if necessary, associate the matching edge $f_i$ with the skeleton edge $e_i$ for each $1\le i\le s$.
Here this association only specifies that $f_i$ will
be used in the construction corresponding to $e_i$.
The skeleton edges $e_1,\ldots,e_s$ need not be pairwise vertex-disjoint,
since they are edges of the same $P_4$ and may share endpoints. In contrast, the edges $f_1,\ldots,f_s$ are pairwise vertex-disjoint because $\mathcal E$ is a matching.
After these assignments,
exactly $3-s=3-|\mathcal E|$ skeleton edges remain unassigned.
For each such edge,
we choose one distinct unused vertex of $C$ for the corresponding construction.

Since $2k$ vertices of $C$ have already been used, there are at least $|C|-2k$ unused vertices available in $C$. By the assumption
$|C|-2k\ge 3-|\mathcal E|$,
the number of available vertices in $C$ is at least the number of remaining skeleton edges. Hence we may assign a distinct unused vertex of $C$ to each of these $3-|\mathcal E|$ skeleton edges.
For a skeleton edge $uv$ with $u\in X$ and $v\in Y$, form its
five-cycle as follows:
\begin{itemize}
\item If the assigned matching edge is $ab$ and both its endpoints
are complete to $Y$, use $uvabyu$ with fresh $y\in Y$.
\item If the assigned matching edge is $ab$ and both its endpoints
are complete to $X$, use $uvxabu$ with fresh $x\in X$.
\item If a vertex $z\in C$ was reserved instead, use
$uvxzyu$ with fresh $x\in X,y\in Y$.
\end{itemize}
If a matching edge satisfies both completeness conditions, use the
first construction. Every edge displayed above exists by the hypotheses.

Every other skeleton edge has one endpoint $c\in C$ and the other
in a partite set. If that endpoint is $v\in Y$, use the five-cycle
$cvxyx'c$
with fresh $x,x'\in X,y\in Y$; interchange $X,Y$ when the endpoint
lies in $X$.

We choose all auxiliary vertices successively,
always avoiding the skeleton vertices, all endpoints
of $\mathcal E$, and all auxiliary vertices
chosen earlier.
At every stage, the number of vertices already
used in the partial embedding is at most $h_5$,
since the completed construction contains exactly
$h_5$ vertices. In addition, the endpoints of
$\mathcal E$ contribute at most six reserved
vertices, some of which may already have been
used in the partial embedding.
Consequently, at most $h_5+6$ vertices
are forbidden in either of $X$ and $Y$ at any step.
Since
\[
|X|,|Y|\ge \beta=2f+10>h_5+6,
\]
we can always choose the required fresh vertices.
Each matching edge is used in exactly one replacing
cycle, and every vertex of $C$ outside the skeleton
is used in at most one replacing cycle.
Thus the internal vertices of the replacing cycles
are pairwise disjoint and are disjoint from all
skeleton vertices. The construction therefore
gives a genuine copy of $H_5$.

In each of the above constructions, the resulting five-cycle
contains a non-skeleton $X$--$Y$ edge. Claim~\ref{claim:3.2}
therefore extends this copy to the prescribed $F^o$, as required.
\end{proof}

\begin{claim}\label{claim:3.4}
The graph $G_1\cup G_2$ has matching number at most two.
\end{claim}

\begin{proof}
Suppose there is a matching $\mathcal E$ of size three in
$G_1\cup G_2$. For $S\subseteq R_0$, set
\[
V_{i,S}=\{v\in V(G_i):N_G(v)\cap R_0=S\}.
\]
The isomorphisms guaranteed by Definition~\ref{def:2.2}(ii)
preserve neighborhoods in $R_0$. Thus every
nonempty $V_{i,S}$ meets each copy, and
$|V_{i,S}|\ge n_i/|H_i|\ge \beta$ for large $n$.

If $V_{1,S}$ and $V_{2,T}$ are nonempty and $|S\cap T|\ge2k$,
choose $C\subseteq S\cap T$ of size $2k$. Every edge of $\mathcal E$
in $G_1$ has both endpoints complete to $V_{2,T}$, and every such
edge in $G_2$ has both endpoints complete to $V_{1,S}$.
The sets $V_{1,S},V_{2,T},C$ are pairwise disjoint, and $C$ is
complete to both classes. Claim~\ref{claim:3.3} would
therefore give $F^o\subseteq G$. Consequently $|S\cap T|\le2k-1$,
and hence
\[
|S|+|T|=|S\cup T|+|S\cap T|\le r_0+2k-1
\]
for every pair of nonempty classes.

For disjoint vertex sets $U,W\subseteq V(G)$,
let $e_G(U,W)$ denote the number of edges between $U$ and $W$.
For $v\in V(G)$ and $U\subseteq V(G)$, let
$d_U(v):=|N_G(v)\cap U|$.
For any $u\in V(G_1)$ and $v\in V(G_2)$, let
\[
S=N_G(u)\cap R_0
\quad\text{and}\quad
T=N_G(v)\cap R_0.
\]
Then $V_{1,S}$ and $V_{2,T}$ are nonempty, since
$u\in V_{1,S}$ and $v\in V_{2,T}$. Hence, by the preceding argument,
$|N_{R_0}(u)\cap N_{R_0}(v)|=|S\cap T|\le 2k-1$.
Since $N_{R_0}(u)\cup N_{R_0}(v)\subseteq R_0$,
we have
\[
d_{R_0}(u)+d_{R_0}(v)
=|N_{R_0}(u)\cup N_{R_0}(v)|
 +|N_{R_0}(u)\cap N_{R_0}(v)|
\le r_0+2k-1.
\]
Summing this inequality over all pairs $(u,v)\in V(G_1)\times V(G_2)$, we obtain
$$
n_2\sum_{u\in V(G_1)}d_{R_0}(u)
+
n_1\sum_{v\in V(G_2)}d_{R_0}(v)
\le
(r_0+2k-1)n_1n_2.
$$
Indeed, each term $d_{R_0}(u)$ with $u\in V(G_1)$ occurs once for each vertex of $G_2$, and hence $n_2$ times, while each term $d_{R_0}(v)$ with $v\in V(G_2)$ occurs $n_1$ times. Dividing by $n_1n_2$, we get $a_1+a_2\le r_0+2k-1$,
where $a_i={n_i}^{-1}\sum_{v\in V(G_i)}d_{R_0}(v)$ for each $i\in\{1,2\}$.

Moreover, each edge between $R_0$ and $G_i$ is counted exactly once in $\sum_{v\in V(G_i)}d_{R_0}(v)$, and hence
$e_G(R_0,G_i)=n_i a_i.$
Therefore
$e_G(R_0,G_1\cup G_2)=n_1a_1+n_2a_2.$
Since $n_i=n/2+O(1)$ and
$0\le a_i\le r_0\le q,$
where $q$ is independent of $n$, we have
$n_1a_1+n_2a_2=\frac n2(a_1+a_2)+O(1).$
Combining this with $a_1+a_2\le r_0+2k-1$, we obtain
$$
e_G(R_0,G_1\cup G_2)
\le
\frac{r_0+2k-1}{2}n+O(1).
$$
Also $n_1n_2=n^2/4-r_0n/2+O(1)$. Together with
\eqref{eq:5A}, these estimates yield
\begin{align*}
e(G)
&=n_1n_2+e_G(R_0,G_1\cup G_2)+e(G_1)+e(G_2)+e(G[R_0])\\
&\le\frac{n^2}{4}+\left(k+\frac12\right)n+O(1),
\end{align*}
which contradicts \eqref{eq:4A}.
\end{proof}

If either $H_i$ had an edge, three of its copies would give a
matching of size three. Thus $H_1=H_2=K_1$. Set
$A_0=V(G_1)$ and $B_0=V(G_2)$, with $a=|A_0|$ and $b=|B_0|$.
These sets are independent and completely joined, with
$a,b=n/2+O(1)$ and $a,b\ge \beta$. By Definition~\ref{def:2.2}(ii), each vertex
of $R_0$ is either complete or anticomplete to each of $A_0,B_0$.
Partition $R_0$ into $C,A_1,B_1,Z$ according as its neighborhood
in $A_0\cup B_0$ is $A_0\cup B_0$, $B_0$, $A_0$, or $\varnothing$,
respectively. Write $c=|C|$ and $z=|Z|$.

\begin{claim}\label{claim:3.5}
We have $c=2k+2$ and $Z=\varnothing$.
\end{claim}

\begin{proof}
We first show that
$c\le 2k+2.$
Indeed, suppose that $c\ge 2k+3$. By definition, every vertex of $C$ is complete to both $A_0$ and $B_0$, while $A_0$ and $B_0$ are completely joined.
Since $|A_0|,|B_0|\ge \beta$, the sets $A_0$ and $B_0$ satisfy the size requirement in Claim~\ref{claim:3.3}.
Taking
$X=A_0$,$Y=B_0$, $\mathcal E=\varnothing$,
and using the fact that $C$ is complete to $A_0\cup B_0$ while $A_0$ is complete to $B_0$, Claim~\ref{claim:3.3} yields a copy of $F^o$ in $G$, contradicting the assumption that $G$ is $F^o$-free. Hence
$c\le 2k+2.$

We next estimate $e(G)$. Write
$a=|A_0|$, $b=|B_0|$.
Since
$V(G)=A_0\mathbin{\dot\cup}B_0\mathbin{\dot\cup}R_0,$
we have
$a+b=n-r_0.$
Moreover, by the structural theorem,
$|a-b|=O(1),$
and $r_0=O(1)$.

By the definitions of $C,A_1,B_1$, and $Z$, all edges between $A_0$ and $B_0$ are present; every vertex of $C$ is complete to both $A_0$ and $B_0$; every vertex of $A_1$ is complete to $B_0$ and anticomplete to $A_0$; and every vertex of $B_1$ is complete to $A_0$ and anticomplete to $B_0$. Therefore
$$
e(G)
=
ab+c(a+b)+|A_1|b+|B_1|a+e(G[R_0]).
$$

We estimate these terms separately.
Since $a+b=n-r_0$, it follows that
$$ab=\frac{(n-r_0)^2-(a-b)^2}{4}=\frac{n^2}{4}-\frac{r_0}{2}n+\frac{r_0^2-(a-b)^2}{4}.$$
As $r_0=O(1)$ and $|a-b|=O(1)$, the last term is $O(1)$. Hence
$ab=\frac{n^2}{4}-\frac{r_0}{2}n+O(1).$

Next,
$c(a+b)=c(n-r_0)=cn-cr_0.$
Since both $c$ and $r_0$ are bounded independently of $n$, we have
$c(a+b)=cn+O(1).$

Similarly, from
$a,b=\frac n2+O(1),$
and the fact that $|A_1|$ and $|B_1|$ are bounded, we obtain
$|A_1|b=\frac{|A_1|}{2}n+O(1)$
and
$|B_1|a=\frac{|B_1|}{2}n+O(1).$
Thus
$|A_1|b+|B_1|a=\frac{|A_1|+|B_1|}{2}n+O(1).$

Finally, since $r_0=O(1)$,
$e(G[R_0])\le \binom{r_0}{2}=O(1).$
Combining the above estimates yields
$$
\begin{aligned}
e(G)
=
\frac{n^2}{4}
-\frac{r_0}{2}n
+cn
+\frac{|A_1|+|B_1|}{2}n
+O(1)
=
\frac{n^2}{4}
+
\left(
-\frac{r_0}{2}
+c
+\frac{|A_1|+|B_1|}{2}
\right)n
+O(1).
\end{aligned}
$$
Since
$r_0=c+|A_1|+|B_1|+z,$
the coefficient of $n$ simplifies as
$$
\begin{aligned}
-\frac{r_0}{2}+c+\frac{|A_1|+|B_1|}{2}
=-\frac{c+|A_1|+|B_1|+z}{2}+c+\frac{|A_1|+|B_1|}{2}
=\frac{c-z}{2}.
\end{aligned}
$$
Therefore,
$e(G)=\frac{n^2}{4}+\frac{c-z}{2}n+O(1).$

On the other hand, by~\eqref{eq:4A}, we have
$$
e(G)
\ge
\frac{n^2}{4}
+(k+1)n
+O(1)
=
\frac{n^2}{4}
+\frac{2k+2}{2}n
+O(1).
$$
Comparing the two estimates gives
$\frac{2k+2-c+z}{2}n\le O(1).$
Equivalently, there exists a constant $M>0$, independent of $n$, such that
$(2k+2-c+z)n\le M.$

Recall that $c\le 2k+2$ and $z\ge0$.
Hence $2k+2-c+z$ is a nonnegative integer.
If it were positive, then it would be at least $1$, and hence
$(2k+2-c+z)n\ge n,$
which is impossible for sufficiently large $n$. Therefore
$2k+2-c+z=0.$
Since both $2k+2-c$ and $z$ are nonnegative, they must both vanish. Thus
$c=2k+2$ and $z=0$.
Hence $Z=\varnothing$, as claimed.
\end{proof}

Set $A=A_0\cup A_1$ and $B=B_0\cup B_1$. Then
$V(G)=C\mathbin{\dot\cup}A\mathbin{\dot\cup}B$ and
$|A|+|B|=n-2k-2$.

\begin{claim}\label{claim:3.6}
The sets $C,A,B$ are independent.
\end{claim}

\begin{proof}
Let $I$ be the independent set given by Lemma~\ref{lem:2.1}(ii), and let $M$ be the graph obtained by cracking all vertices of $I$ and choosing Type~I for every edge in $\delta_F(I)$.
By Proposition~\ref{prop:2.4}, we have
$M\in\mathcal M(F^o)$.

We first determine the structure of $M$.
By Lemma \ref{lem:2.1}(ii), we have
$F-I\cong P_4\cup(2k-2)K_1,$
where the unique nontrivial component is the path
$P=x_{2k}x_{2k+1}y_{2k+1}y_1.$
Since $F$ is cubic, under the Type~I operation each vertex $v\in V(F-I)$ acquires exactly
$3-d_{F-I}(v)$
new pendant neighbors. Consequently, each isolated vertex of $F-I$ becomes the center of a copy of $K_{1,3}$.
Hence
$M=Q\cup(2k-2)K_{1,3},$
where $Q$ is the tree obtained from the path $x_{2k}x_{2k+1}y_{2k+1}y_1$ by attaching two leaves to each endpoint and one leaf to each internal vertex.
The graph $Q$ is depicted in Figure~\ref{fig:Q-typeI},
where its six pendant vertices are labelled $a_1,a_2,b,c,d_1,d_2$.

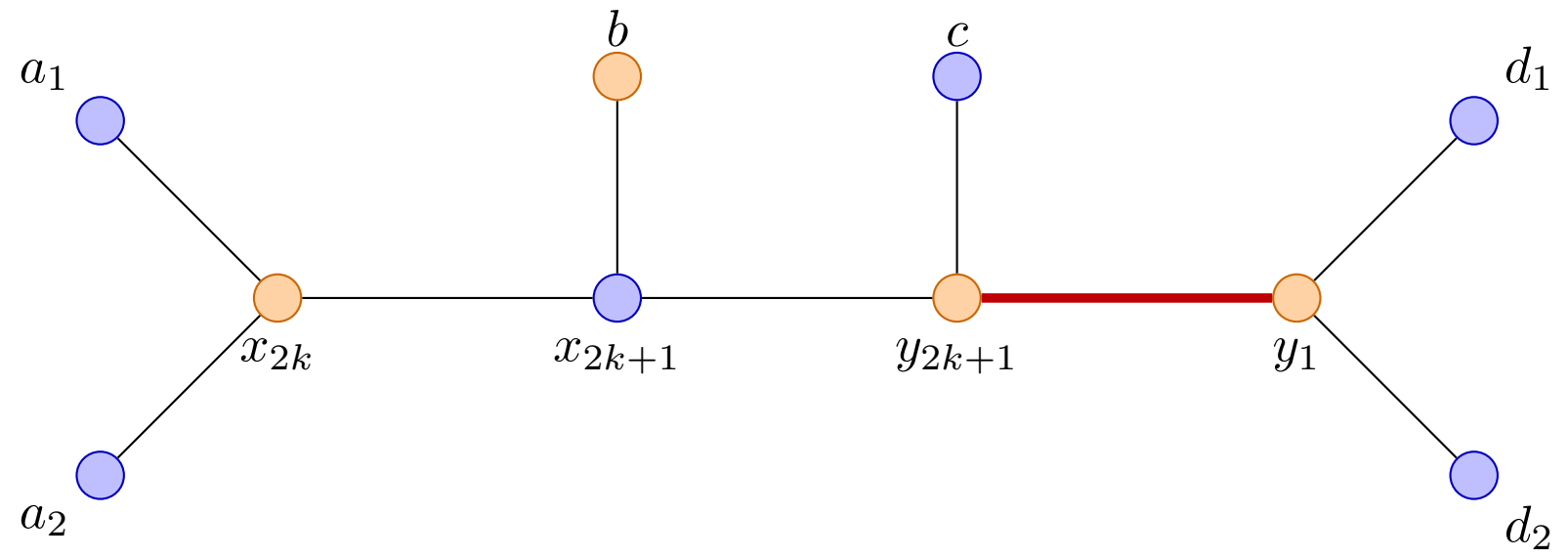
\begin{figure}[htbp]
\centering
\begin{tikzpicture}[
    scale=1,
    cvertex/.style={
        circle,
        draw=orange!80!black,
        fill=orange!35,
        minimum size=3.2mm,
        inner sep=0pt
    },
    xvertex/.style={
        circle,
        draw=blue!75!black,
        fill=blue!25,
        minimum size=3.2mm,
        inner sep=0pt
    },
    every edge/.style={draw=black, line width=0.9pt},
    font=\small
]

% Core path vertices
\node[cvertex] (x2k) at (0,0) {};
\node[xvertex] (x2kp1) at (2.3,0) {};
\node[cvertex] (y2kp1) at (4.6,0) {};
\node[cvertex] (y1) at (6.9,0) {};

% Core path labels
\node[below=4pt] at (x2k) {$x_{2k}$};
\node[below=4pt] at (x2kp1) {$x_{2k+1}$};
\node[below=4pt] at (y2kp1) {$y_{2k+1}$};
\node[below=4pt] at (y1) {$y_1$};

% Core path edges
\draw (x2k) -- (x2kp1);
\draw (x2kp1) -- (y2kp1);
\draw[red!75!black, line width=1.8pt] (y2kp1) -- (y1);

% Leaves at x_{2k}
\node[xvertex] (a1) at (-1.2,1.2) {};
\node[xvertex] (a2) at (-1.2,-1.2) {};
\draw (x2k) -- (a1);
\draw (x2k) -- (a2);
\node[above left=2pt] at (a1) {$a_1$};
\node[below left=2pt] at (a2) {$a_2$};

% Leaf at x_{2k+1}
\node[cvertex] (b) at (2.3,1.5) {};
\draw (x2kp1) -- (b);
\node[above=2pt] at (b) {$b$};

% Leaf at y_{2k+1}
\node[xvertex] (c) at (4.6,1.5) {};
\draw (y2kp1) -- (c);
\node[above=2pt] at (c) {$c$};

% Leaves at y_1
\node[xvertex] (d1) at (8.1,1.2) {};
\node[xvertex] (d2) at (8.1,-1.2) {};
\draw (y1) -- (d1);
\draw (y1) -- (d2);
\node[above right=2pt] at (d1) {$d_1$};
\node[below right=2pt] at (d2) {$d_2$};

\end{tikzpicture}
\caption{The graph $Q$.}
\label{fig:Q-typeI}
\end{figure}

Set $X=A_0$ and $Y=B_0$.
Then $X$ and $Y$ are independent and completely joined,
$|X|,|Y|\ge \beta$, and $C$ is complete to $X\cup Y$.
As illustrated in Figure~\ref{fig:Q-typeI}, we place $x_{2k},y_{2k+1},y_1,b$ in $M_C$ and $x_{2k+1},a_1,a_2,c,d_1,d_2$ in $M_X$. In the figure, orange vertices represent the vertices of $M_C$, blue vertices represent those of $M_X$,
and the red edge is the unique edge contained in $M[M_C]$.
For each of the remaining $2k-2$ copies of $K_{1,3}$ in $M$, place its center in $M_C$ and its three leaves in $M_X$. Then $|M_C|=4+(2k-2)=2k+2$ and $|M_X|=6+3(2k-2)=6k$.
Moreover, $M_X$ is independent, $M[M_C]\cong K_2\cup 2kK_1$,
and every edge of $M$ other than $y_{2k+1}y_1$
has one endpoint in $M_C$ and the other in $M_X$.

Now suppose that $c_1c_2\in E(G[C])$. Since $|C|=|M_C|=2k+2$, we may embed $M_C$ bijectively into $C$, sending the unique edge $y_{2k+1}y_1$ to $c_1c_2$. Also, since $|M_X|=6k<\beta\le |X|$, we may embed $M_X$ injectively into $X$. As $C$ is complete to $X$, every edge of $M$ between $M_C$ and $M_X$ is present under this embedding.
Hence $M\subseteq G[C\cup X].$

Recall that $M\in\mathcal M(F^o)$. By the definition of the decomposition family and the observation established above, we have
$F^o\subseteq (M\cup E_f)\vee E_f$.
We now show that $G$ contains $(M\cup E_f)\vee E_f$.
The embedding of $M$ constructed above uses exactly $6k$ vertices of $X$.
Since $f\ge h_5=22k+11>6k$ and $|X|,|Y|\ge \beta=2f+10$,
we can choose an independent set $X'\subseteq X\setminus V(M)$ of size $f$ and an independent set $Y'\subseteq Y$ of size $f$.
Recall that both $X$ and $Y$ are independent and are completely joined.
Moreover, $Y$ is complete to $C\cup X$, and hence to $V(M)\cup X'$.
Consequently, the vertices of $X'$ and $Y'$ together with the embedded copy of $M$ form a copy of $(M\cup E_f)\vee E_f$ in $G$.
It follows that $F^o\subseteq G$, contradicting the assumption that $G$ is $F^o$-free. Therefore $G[C]$ is independent.

Next, $A_0$ is independent and $A_1$ is anticomplete to $A_0$,
so every edge of $G[A]$ would lie in $G[A_1]$. If such an edge
$ab$ existed, take $\mathcal E=\{ab\}$. Its endpoints are both
complete to $B_0$, and $|C|+|\mathcal E|=2k+3$.
Claim~\ref{claim:3.3}, with $X=A_0$ and $Y=B_0$, would
therefore give $F^o\subseteq G$, a contradiction. Thus $G[A]$ is
independent. Interchanging $A$ and $B$ proves that $G[B]$ is
independent as well.
\end{proof}

We can now complete the proof by a final edge count.
Since $C,A,B$ are all independent, it follows that
\begin{align*}
e(G)
\le (2k+2)(|A|+|B|)+|A|\cdot |B|
\le (2k+2)(n-2k-2)+\Big\lfloor\frac{(n-2k-2)^2}{4}\Big\rfloor
 =g_k(n).
\end{align*}
Together with the lower bound~\eqref{eq:4A}, equality holds throughout.
Since the lower bound is attained by $E_{2k+2}\vee T_{n-2k-2,2}$ and $G$ is extremal,
neither of these inequalities can be strict. In the first line, the
deficit is exactly the number of missing edges between distinct parts
$C,A,B$, so every such edge is present. For the second line, use
$4|A||B|=(|A|+|B|)^2-(|A|-|B|)^2$ with $|A|+|B|=n-2k-2$; the integral product is
maximal exactly when $\big||A|-|B|\big|\leq 1$. Therefore
$G\cong E_{2k+2}\vee T_{n-2k-2,2}$.
Since $G$ was arbitrary in $\mathrm{EX}(n,F^o)\cap\mathcal D(n,2,q)$,
the uniqueness assertion of Lemma~\ref{lem:2.6} now shows
that this is the unique extremal graph for all sufficiently large $n$.
\end{proof}

\end{document}